\documentclass[a4paper,10pt]{amsart}
\usepackage{amsfonts}
\usepackage{amsmath}
\usepackage{amssymb}
\usepackage{amsthm}
\title{Explicit Factorization of  Prime Integers in Quartic Number Fields defined by $X^4+aX+b$}
\author{L houssain El Fadil}
\date{}
\thanks{Supported by ERCIM}
\newtheorem{teor}{Theorem}[section]

\newtheorem{lem}[teor]{Lemma}

\newtheorem{exs}[teor]{Examples}

\newcommand{\z}{\mathbb Z}
\newcommand{\q}{{\mathbb Q}}

\newcommand{\F}{\mathbb F}
\def\op{\operatorname}

\def\al{\alpha}
\def\as#1{\renewcommand\arraystretch{#1}}

\def\be{\bigskip}

\def\diso{\lower.4ex\hbox{$\downarrow$}\raise.4ex\hbox{\mbox{\scriptsize $\wr$}}}

\def\fp{\mathbb F_p}
\def\fph{\mathbb{F}_{\ph}}

\def\al{\alpha}

\def\iso{\,\lower .6ex\hbox{$\stackrel{\lra}{\mbox{\tiny $\sim\,$}}$}\,}

\def\la{\lambda}

\def\lg{l\raise.6ex\hbox to.2em{\hss.\hss}l}
\def\lra{\longrightarrow}

\def\md#1{\ \mbox{\rm(mod }{#1})}

\def\nph#1{N_{\ph}(#1)}

\def\orb{\hbox to  .3em{$\backslash$}\backslash}

\def\p{\mathfrak{p}}

\def\ph{\phi}

\def\K{{\mathbb{K}}}

\def\rd{\op{red}}

\def\t{\theta}

\newcounter{cs}
\stepcounter{cs}
\newcommand{\casos}{\begin{itemize}}
\newcommand{\fcasos}{\end{itemize}\setcounter{cs}{1}}

\newfont{\tit}{cmr12 scaled \magstep3}
\begin{document}
\maketitle
\begin{abstract}
For every prime integer $p$, an explicit factorization   of the
principal ideal $p\z_K$ into prime ideals of $\z_K$ is given,
where $K$ is a quartic number  field defined by an irreducible
polynomial $X^4+aX+b\in\z[X]$.
\end{abstract}
\section*{Introduction}
 Let $K$ be a quartic number  field  defined by
an irreducible polynomial $P(X)=X^4+aX+b\in\z[X]$, $\z_K$ its ring
of integers, $\triangle$ the discriminant of $P$, $d_K$ its
discriminant, $\al$ a complex root of $P(X)$ and
$ind(P)=[\z_K:\z[\alpha]]$ the index of $\z[\alpha]$. In this
paper, the goal is  to give an explicit factorization  of $p\z_K$
into prime ideals of $\z_K$; the form
$p\z_K=\prod_{i=1}^rP_i^{e_i}$ and for every prime  factor $P_i$,
an integral element $w_i$ such that $P_i=(p,w_i)$ are given. Let
$p$ be a prime integer. It is well known that if $p$ does not
divide $ind(P)$, then  the Dedekind's theorem gives us an explicit
factorization of the principal ideal $p\z_K$ into prime ideals of
$\z_K$:
  $p\z_K$ is $p$-analogous to the factorization  of $\bar P(X)$ modulo $p$.
(see, for example \cite [page 257]{Co}).
\begin{teor}
 Let $p$ be a  prime integer. Denote by $\bar {}$ the canonical map of $\z[X]$ into ${F_p}[X]$, and let $\bar  P(X)=\prod_{i=1}^rg_i(X)^{e_i}$,
 where $g_1(X)$,..,$g_r(X)$ are distinct irreducible in ${F_p}[X]$
 and $e_1$,..$e_r$ are positive integers.
 For every $i$, let $P_i=(p,f_i(\al))$, where $f_i\in\z[X]$ is a monic lifting over $g_i(X)$. Then\\
If $p$ does not divides $ind(P)$, then
$p\z_K=\prod_{i=1}^rP_i^{e_i}$ and for every $i$, $e(P_i/p)=e_i$ and
$f(P_i/p)=f_i=deg(g_i)$.
\end{teor}

If $p$ is not a common index divisor of $K$, then there exist an
element $\phi\in\z_K$ which generates $\z_K$ and
$v_p(ind(\pi_{\phi})=0$, where $\pi_{\phi}$ is the minimal
polynomial of $\phi$, and then we can apply Dedekind's theorem to
obtain the prime ideal decomposition explicitly.  However given
$\al$, it is not easy to determine such an element $\phi$ in
general. The construction of $\phi$ was based on the $p$-integral
bases of  $\z_K$ given in \cite{AW}. If $p$  is common index
divisor of $K$, then for every prime $P$ factor of  $p\z_K$,  an
element  $\beta\in \z_K$ such that $v_P(\beta)=1$ and for every
prime ideal $Q\neq P$ above $p$, $v_Q(\beta)=0$ will be
constructed. A such  $\beta\in \z_K$ satisfies  $P=(p, \beta)$.

For every prime $p$, let $v_p$ be the $p$-adic discrete valuation
defined in $\q_p$ by $v_p(p)=1$. $v_p$ is extended to $\q_p[X]$ by
$v_p(A(X)):=Min\{v_p(a_i)\, ,\, 0\le i\le r\}$, where
$A(X)=\sum_{i=0}^ra_iX^i$. For every $(x,m)\in\z^2$, denote
$x_p=\frac{x}{p^{v_p(x)}}$, $x\, \md{m}$ the remainder of the
Euclidean division of $x$ by $m$. For every odd prime $p$, denote $(\frac{-}{p})$ the Legendre
symbol.  If $p$ does not divide $x$ and there exists $t\in\z$ such
that $t^n=x\md{p}$, then let $(\frac{x}{p})_n=1$. Otherwise,
$(\frac{x}{p})_n\neq 1$.
\section{Newton polygons}
Let $F(X)$ be an irreducible polynomial in $\z[X]$, $\phi\in
\z[X]$ a monic polynomial of degree at least $1$. Let
$F(X)=\sum_{i=0}^na_i(X)\phi(X)^{i}=
a_0(X)+..+a_{n-1}(X)\phi(X)^{n-1}+a_{n}(X)\phi(X)^n$ be  the
$\phi$-adic development of $F(X)$ ( for every $i$, $deg(a_i(X))\le
deg(\phi)-1$). Let $p$ be a prime integer. The $\phi$-Newton
polygon of $F(x)$, with respect to $p$, is the lower convex
envelope of the set of points $(i,u_i)$, $u_i<\infty$, in the
Euclidian plane, where $u_i=v_p(a_i(X))$.
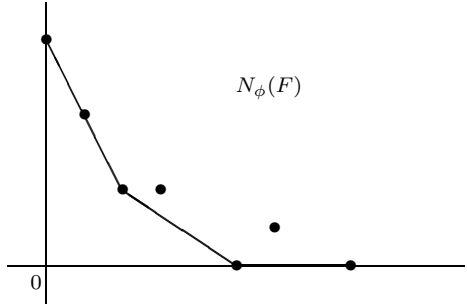
\begin{center}
\setlength{\unitlength}{5.mm}
\begin{picture}(12,8)
\put(7.85,-.15){$\bullet$}\put(5.85,.85){$\bullet$}\put(4.85,-.15){$\bullet$}\put(2.85,1.85){$\bullet$}
\put(1.85,1.85){$\bullet$}\put(.85,3.85){$\bullet$}\put(-.15,5.85){$\bullet$}
\put(-1,0){\line(1,0){12}}\put(0,-1){\line(0,1){8}}
\put(5,0){\line(-3,2){3}}\put(2,2){\line(-1,2){2}}\put(5,.03){\line(-3,2){3}}
\put(2,2.03){\line(-1,2){2}}\put(8,0){\line(-1,0){3}}\put(8,.02){\line(-1,0){3}}
\put(4,-.8){\begin{footnotesize}\end{footnotesize}}
\put(-.4,-.6){\begin{footnotesize}$0$\end{footnotesize}}
\put(5,4.6){\begin{footnotesize}$N_{\phi}(F)$\end{footnotesize}}
\end{picture}\qquad\qquad
\end{center}
\be
The $\phi$-Newton polygon is the union of different adjacent sides $S_0,\dots,S_g$ with
increasing slope $\la_0<\la_1<\cdots<\la_g$.
We shall write $N_{\phi}(F)=S_0+\cdots+S_g$. The principal part of $N_{\phi}(F)$,
 denoted $N^{+}_{\phi}(F)$, is the polygon determined by the sides of negative slope of
 $N_{\phi}(F)$.

For every $0\le i\le n$, we attach to any abscissa  the following {residual coefficient}
$c_i\in\fph$:
$$\as{1.6}
c_i=\left\{\begin{array}{ll} 0,&\mbox{ if $(i,u_i)$ lies strictly
above $N$ or
}u_i=\infty,\\
\rd\left(\dfrac{a_i(X)}{p^{u_i}}\right),&\mbox{ if
$(i,u_i)$ lies on }N.
\end{array}
\right.
$$
Let $S$ be one of the sides of $N$, with slope $\la$, and let
$\la=-h/e$, with $h,e$ positive coprime integers. Let $l=\ell(S)$
be the length of the projection of $S$ to the $x$-axis,
$d(S):=\ell(S)/e$: the degree of $S$.  Note that $S$ is divided
into $d(S)$ segments by the points of integer coordinates that lie
on  $S$. Let $s$ be the initial abscissa of $S$ and $d$ the degree
of $S$. If $\phi$ is a monic polynomial such that $\bar\phi$ is
irreducible factor of $\bar F(X)$ modulo $p$, then let
$\fph:=\frac{\fp[X]}{(\phi)}$ and
$F_S(Y):=c_s+c_{s+e}\,Y+\cdots+c_{s+(d-1)e}\,Y^{d-1}+c_{s+de}\,Y^d\in\fph[Y]$,
the residual polynomial of $f(X)$ attached to $S$.\\

In the remainder, $F(X)$ is a monic  irreducible polynomial, $\t$
is a complex root of $F(X)$ and $\phi$ is a monic polynomial such
that $\bar\phi$ is irreducible factor of $\bar F(X)$ modulo $p$.
Let $N=S_k+..+S_1+S_0$ be the $X\phi$-Newton polygon of $F(X)$
with respective slopes $\lambda_k<..<\lambda_1<\lambda_0$. $F(X)$
is said to be $\phi$-regular if for every $i$, $F_{S_i}(Y)$ is
square free.
\section{Main results}

 \be In this section, for every prime integer $p$, an explicit factorization   of
the principal ideal $p\z_K$ into prime ideals of $\z_K$ is given
(Note only the form is given, but two element generators of each
prime ideal factor are given). Recall the  celebrated Theorem of
Hensel:
\begin{teor}[Hensel]\label{hensel}
Let $P(X)\in\z[X]$ be a monic irreducible polynomial, $p$  a prime
integer, $\q_p$ the $p$-adic completion of $\q$ and $\z_p$ its
ring of integers. Let $P(X)=F_1(X)\cdots F_t(X)$ be the
factorization into a product of monic irreducible polynomials in
$\q_p[X]$ and consider the local fields $K_i=\q_p[X]/(F_i(X))$,
for $i=1,\dots,t$. The factorization of $p\z_K$ into a product of
prime ideals of $\z_K$ is: $p\z_K=\p_1^{e_1}\cdots \p_t^{e_t}$,
where $f(\p_i)=f(K_i/\q_p)$ and $e_i=e(K_i/\q_p).$
\end{teor}
\begin{lem}
Assume that  $P(X)=\bar \phi^l(X)\,\md{p}$ such that $N_{\phi}{P}=S$
is one side  and $P_S(Y)$ is irreducible in $\fph[Y]$, where
$\fph=\frac{\fp[X]}{(\phi(X))}$. Let $\la=-h/e$ be the  slope of
$S$ such that $h$ and $e$ are positive coprime integers. Define
$\al_S:=\ph(\al)^{e}/p^h$. Then $v_p(\al_S)=0$ and $P_S(\al_S)=0$.
\end{lem}
\begin{proof}
Let $P(X)=\sum_{k=0}^la_k(X)\phi(X)^k$ be the $\phi$-adic
development of $P(X)$. Since $\nph{P}$ is one side of slope
$-\la$, we have  for every $i$, $v(a_{i}(X))\ge (i-l)\la$ and
$v(a_0(X))=-l\la$. Thus, $v_p(\ph(\al))=-\la$. Moreover, since for
every $0\le j\le d$ and for every $je<k<(j+1)e$,  $\rd
\frac{a_k(X)}{p^{u_k}}=0$ and $\frac{P(\al)}{p^l}=0$, then
$\overline{
(\frac{\phi(\al)^e}{p^h})^d+\frac{a_{l-e}(\al)}{p^h}\frac{\phi(\al)^e}{p^h})^{d-1}+\cdots
+\frac{a_{e}(\al)}{p^{h(d-1)}}\frac{\phi(\al)^e}{p^h})+\frac{a_{le}(\al)}{p^{l}}}=0$
in  $\fph$. Therefore, $P_S(\al_S))=0$ in  $\fph$.
\end{proof}

 Let $P(X)\in\z[X]$ be a monic irreducible polynomial and
  $\bar P(X)=\prod_{i=1}^r\phi_i^{l_i}\,\md{p}$ its factorization into irreducible polynomials of $\fp[X]$. For every
$i:=1..r$, let $N_i:=N^+_{\phi_i}(P)=S_1^i+\cdots+S_{k_i}^i$ and
for every $j:=1..k_i$, let
$P_{S_j}(Y)=\prod_{s=1}^{r_{ij}}\psi_s(Y)^{n_s}$ be the
factorization of $P_{S_j}(Y)$ into irreducible polynomials of
$\fph[Y]$, where $\fph=\frac{\fp[X]}{(\phi(X))}$.
\begin{teor}
 Under these hypothesis, if every $P_{S_j}(Y)$ is
square free (every $n_s=1$), then
$$p\z_K=\prod_{i=1}^{r}\prod_{j=1}^{k_i}\prod_{s=1}^{r_{ij}}\p_s^{e(S_j^i)},$$
 and  for every $i:=1..r$, $j:=1..k_i$ and  $s:=1..r_{ij}$,
 $v_{\p_s}({\phi_i}(\al))=-e(S_j^i)\la_j$, where $e(S_j^i)$ is the ramification index of the
 side $S_j^i$ and $-\la_j$ is its slope.
\end{teor}
\begin{proof}
By Hensel Theorem, it suffices to factorize $P(X)$ into
irreducible polynomials of $\q[X]$. By Hensel Lemma, it suffices
to show this result for $r=1$. By  Theorem of the polygon we can
assume that $N_1=S$ is one side. Since $P_S(Y)$ is square free, by
Theorem of the residual polynomial, we can assume that $P_S(Y)$ is
irreducible in $\fph[Y]$. By Theorem of the product, if $P_S(Y)$
is irreducible in $\fph[Y]$, then $P(X)$ is irreducible in
$\q[X]$.

Now, assume that $P(X)\in \z_p[X]$ is a monic irreducible
polynomial
 in $\q[X]$ such that  $\bar P(X)=\phi^l(X)\,\md{p}$,
$N_{\phi}(P)=S$ is one side and $P_S(Y)$ is irreducible in
$\fph[Y]$ of degree $d(S)$. Denote $\K=\q_p[\al]$, $\z_{\K}$ its
ring of integers over $\z_p$ and $\p$ the maximal ideal of
$\z_{\K}$. By Hensel Theorem, $p\z_{\K}=\p^{e(\p)}$, where
$e(\p)=\frac{deg(P)}{f(\p)}$.

Let $\lambda =-\frac{e}{h}$ be the slope of $S$ such $e$ and $h$
are positive coprime. Let $\al_S=\frac{\phi^e(\al)}{p^h}$ and  $
\iota\colon \z_p[X,Y]\lra \z_p,\qquad a(X,Y)\mapsto a(\al,\al_S).$
Then the ideal $\mathfrak{m}=\iota^{-1}(\p)$ is a maximal ideal of
$\z_p[X,Y]$ generated by $(p,\ph(X),\varphi(X,Y))$, where
$\varphi(X,Y)\in \z_p[X,Y]$ is a monic polynomial such that
$\psi(Y):=\rd(\varphi(X,Y))$ is irreducible in $\fph[Y]$. Thus,
$\frac{\z_{\K}}{\p}\simeq \frac{\fph[Y]}{(\psi(Y))}$ and
$f(\p)=d(S).m$. Thus,
$e(\p)=\frac{deg(P)}{d(S).m}=\frac{l.m}{d(S).m}=e(S)$, where
$m=deg(\phi)$.\\
On the other hand, since $N_1=S$ is one side of slope $\la$, we
have $v_{p}(\ph(\al)^l)=v_{p}(a_0(\al))=-l\la$,
$v_{p}(\ph(\al))=-\la$ and $v_{\p}(\ph(\al))=-e(\p)\la$.
\end{proof}

 Let $P(X)=X^4+aX+b$. Note that if there exists a prime $p$ such that $v_p(a)\ge 3$ and $v_p(b)\ge 4$, then let $q_1$ and $q_2$ be respectively the quotient of  $v_p(a)$ by $3$ and of $v_p(b)$ by $4$. Let $\t:=\frac{\al}{p^q}$, where $q:=Min(q_1,q_2)$. Then $\t$ is integral with minimal polynomial $F(X)=X^4+AX+B\in\z[X]$, where $A=\frac{a}{p^{3q}}$ and $B=\frac{b}{p^{4q}}$. As  $K=\q[\t]=\q[\al]$, then up to replace $\al$ by $\frac{\al}{p^q}$, we can assume that  for every prime $p$, $v_p(a)\le 2$ or $v_p(b)\le 3$.
\begin{teor}
Let $p$ be a prime integer.  In the following tables, the form of $p\z_K$
as a product of prime ideals of $\z_K$,   and for every prime
factor $P$ of $p\z_K$  an integral element $\phi$ of $K$ such that
$P=(p,\phi)$ are given:
 If  $p\z_K=\prod_{i}P_i^{e_i}$, then for every $ i\neq j$, an element $\beta_i\in K$
 such that $v_{P_i}=(\beta_i)=1$ and $v_{P_j}=(\beta_i)=0$ is given.
 $$Table A: v_p(a)\ge 1 \, and  \,v_p(b)\ge 1$$
 {\small $$\begin{tabular}{|c|c|c|c|c|c|}
\hline
Case&Conditions& $p$  & $p\z_K$& Generators\\
\hline A1&$v_p(b)=3$, $v_p(a)\ge 3$& & $P^4$&$P=(p,\frac{\al^3}{p^2})$\\
\hline A2&$v_p(b)\ge 3$, $v_p(a)=2$& & $P_1P_2^3$&
$P_1=(p,\frac{\al^3}{p^2})$, $P_2=(p,a_p+\frac{\al^3}{p^2})$\\
\hline A3&$v_p(b)\ge 2$, $v_p(a)=1$& & $P_1P_2^3$&
$P_1=(p,\frac{\al^3}{p})$, $P_2=(p,a_p+\frac{\al^3}{p})$\\
\hline
A4&$v_p(b)= 2$, $v_p(a)\ge 2$& $\neq 2$&  $P^2$& $P=(p,\frac{\al^3}{p})$\\
  &$(\frac{-b_p}{p})=-1$ & &&\\
 \hline
 A5&$v_p(b)= 2$, $v_p(a)= 2$& $\neq 2$&  $P_1^2P_2^2$& $P_1=(p,t+\frac{\al^2}{p})$\\
  &$(\frac{-b_p}{p})=1$ & &&$P_2=(p,-t+\frac{\al^2}{p})$, ($v_p(t^2+b_p)=1$)\\
 \hline
A6&$v_p(b)= 2$, $v_p(a)\ge 3$& $\neq 2$&  $P_1^2P_2^2$& $P_2=(p,t+\frac{\al^3+\al^2}{p})$\\
 &$(\frac{-b_p}{p})=1$ & &&$P_2=(p,-t+\frac{\al^3+\al^2}{p})$, ($t^2+b_p=0\,\md{p}$)\\
 \hline
 A7&$v_p(b)=1$, $v_p(a)\ge 1$& & $P^4$&$P=(p,\al)$\\
\hline
A8&$v_p(b)=2$, $v_p(a)\ge 2$&$2$ & &go to $Table A8$\\
\hline
\end{tabular}$$}
{\small $$Table A8:\, v_2(b)= 2\, and \,v_2(a)\ge 2$$
$$\begin{tabular}{|r|c|c|c|c|r|} \hline
 Case&conditions& $\phi_2$&$2\z_K$&Generators\\
\hline
 A8.1&$v_2(b)= 2$, $v_2(a)= 2$& $\frac{\al^2+2}{2}$& $P^4$&
$P=(2,\phi_2)$\\
\hline
A8.2& $b= 12\md{16}$, $v_2(a)=3$& $\frac{\al^3+2\al}{4}$&
$P^4$&$P=(2,\phi_2)$\\
\hline
A8.3& $b= 4\md{16}$, $v_2(a)=3$& $\frac{\al^3+2\al^2+2\al}{4}$&
$P^4$&
$P=(2,\phi_2)$\\
\hline
A8.4&$b=4\md{32}, v_2(a)=4$ & $\frac{\al^2+2\al+2}{4}$&$P^4$&$P=(2,\phi_2)$\\
\hline
A8.5&$b=20\md{64}, v_2(a)=4$ & $\frac{\al^2+2\al-2}{4}$&$P^2$&$P=(2,\phi_2)$\\
\hline
A8.6&$b=52\md{64}, a=16\md{64}$ & $\frac{\al^3-2\al^2-2\al}{8}$&$P_1^2P_2^2$&$P_1=(2,\phi_2)$, $P_2=(2,\phi_2+1)$\\
\hline
A8.7&$b=52\md{64}, a=48\md{64}$ & $\frac{\al^3-2\al^2+6\al}{8}$&$P_1^2P_2^2$&$P_1=(2,\phi_2)$, $P_2=(2,\phi_2+1)$\\
\hline
A8.8&$b=12\md{32},  v_2(a)\ge 4$ & $\frac{\al^3+2\al}{4}$&$P^2$&$P=(2,\phi_2)$\\
\hline
A8.9&$b=28\md{32},  v_2(a)= 4$ & $\frac{\al^2+2}{4}$&$P_1^2P_2^2$&$P_1=(2,\phi_2)$, $P_2=(2,\phi_2+1)$\\
\hline
A8.10&$b=28\md{32},  v_2(a)\ge 5$ & $\frac{\al^2+12\al+2}{4}$&$P_1^2P_2^2$&$P_1=(2,\phi_2)$, $P_2=(2,\phi_2+1)$\\
\hline
A8.11&$b=20\md{32}, v_2(a)\ge 5$ &$\frac{\al^2+2\al-2}{4}$&$P^4$&$P=(2,\phi_2)$\\
\hline
A8.12&$b=36\md{64},  v_2(a)\ge 5$ &$\frac{\al^2+2\al+2}{4}$&$P^2$&$P=(2,\phi_2)$\\
\hline
A8.13&$b=4\md{64},  v_2(a)= 5$ &$\frac{\al^3-2\al-4}{8}$&$P_1^2P_2^2$&$P_1=(2,\phi_2)$, $P_2=(2,\phi_2+1)$\\
\hline
A8.14&$b=4\md{64},  v_2(a)\ge 6$ &$\frac{\al^3+4\al^2-2\al+4}{8}$&$P_1^2P_2^2$&$P_1=(2,\phi_2)$, $P_2=(2,\phi_2+1)$\\
\hline
\end{tabular}$$}
 {\small $$Table B:\, v_p(b)\ge1\, and \, v_p(a)=0$$
$$\begin{tabular}{|c|c|c|c|c|c|c|}\hline
Case &Conditions&$p$&$\phi$& $p\z_K$& Generators\\\hline
   B1&&$2$ & $\al$& $P_1P_2P_3$& $P_1=(2,\al)$, $P_2=(2,1+\al)$\\
&  &  &  &$P_3=(2,1+\al+\al^2)$\\
\hline
 B2&$(\frac{-a}{p})_3\neq 1$&$\ge
5$ & $\al$& $P_1P_2$& $P_1=(p,\al)$, $P_2=(p,a+\al^3)$\\
\hline
B3&$(\frac{-3}{p})=-1$, $(\frac{-a}{p})_3=1$ &$\ge 5$ & $\al$& $P_1P_2P_3$& $P_1=(p,\al)$, $P_2=(p,-u+\al)$\\
& & &  & & $P_3=(p,u^2+u\al+\al^2)$ \\
&  &&  & & $u^3=-a\,\md{p}$\\
\hline
B4&$(\frac{-3}{p})=1$, $(\frac{-a}{p})_3=1$&$\ge 5$ & $\al$& $P_1P_2P_3P_4$& $P_1=(p,\al)$, $P_2=(p,-u+\al)$\\
  &  && & & $u^3=-a\,\md{p}$,  $P_3=(p,v_1+\al)$\\
 &  && & &    $P_4=(p,-(u+v_1)+\al)$\\
  &  & & & & $v^2=-3\,\md{p}, 2v_1=-u(1+v)\,\md{p}$ \\
  \hline
B5&$ v_3(b)\ge 2, v_3(a)=0$&3 &$\al$&$P_1P_2^3$&$P_1=(3,\al)$, $P_2=(3,\al-a)$\\
&$a^2\neq 1\md{9}$&  &  & &\\
 \hline
 B6&$ v_3(b)\ge 2$, $a^2= 1\md{9}$&3 &&$P_1P_2^2P_3$&go to $Table B6$\\
 \hline
B7& $v_3(a)=0, b=6\md{9}$ &3&$\al$&$P_1^3P_2$&$P_1=(3,\al-a)$, $P_2=(3,\al)$\\
&  $a^2\neq 4\md{9}$ & & & &\\
 \hline
 B8&$ b=6\md{9}, a^2=4\md{9}$ &3  &$\al$&$P_1P_2^2P_3$& go to $Table B6$\\
 \hline
 B9&$v_3(a)=0, b=3\md{9}$ &3  &$\al$&$P_1^3P_2$&$P_1=(3,\al-a)$, $P_2=(3,\al)$\\
 &  $ a^2\neq 7\md{9}$ & & & &\\
  \hline
   B10&$b=3\md{9},  a^2= 7\md{9}$ &3  &&$P_1P_2^3$& $P_1=(3,\t-4a)$, $P_2=(3,\frac{\t^3-4a\t^2}{3})$\\
 &  $v_3(b+a^4-a^2)=2$&  &  & &$\t=\al-a$\\
  \hline
  B11&$b=3\md{9},  a^2= 7\md{9}$ &3  &&& go to $Table B11$\\
  & $v_3(b+a^4-a^2)\ge 3$&  &  & &\\
  \hline
  \end{tabular}$$}
  {\small$$Table B11$$
Let $s\in \z$ such that $as=-4b_3\md{3^{v_3(\triangle)+1}}$, and
let $\theta=\al-s$, $A=4s^3+a$ and $B=s^4+as+b$.
$$\begin{tabular}{|c|c|c|c|c|c|c|}\hline
Case&Conditions&$\phi$&  $3\z_K$& Generators\\\hline
 B11.1&$ v_3(\triangle)=6$, $B_3=1\,\md{3}$ &$\frac{\theta^3+4s\theta^2+6s^2\theta+A}{9}$
&$P_1P_2$&$P_1=(3,\phi)$, $P_2=(3,\phi^3-\phi-s)$\\
 \hline
B11.2&$ v_3(\triangle)=6$ &$\frac{\theta^3+4s\theta^2+6s^2\theta+A}{9}$&$P_1P_2P_3$&$P_1=(3,\phi)$, $P_2=(3,\phi^2-\phi-1)$\\
  &$B_3=-1\,\md{3}$, $s=1\,\md{3}$ && &$P_3=(3,\phi+1)$\\
 \hline
B11.3&$ v_3(\triangle)=6$ &$\frac{\theta^3+4s\theta^2+6s^2\theta+A}{9}$&$P_1P_2P_3$&$P_1=(3,\phi)$, $P_2=(3,\phi^2+\phi-1)$\\
  &$B_3=-1\,\md{3}$, $s=-1\,\md{3}$ && &$P_3=(3,\phi-1)$\\
 \hline
 B11.4&otherwise && go to $Table B11.4$&\\
 \hline
\end{tabular}$$}
{\small  $$Table C: v_p(b)=0 , v_p(a)\ge 1$$
$$\begin{tabular}{|c|c|c|c|c|c|c|c|}
\hline
&Conditions&$p$&$\phi$& $p\z_K$& Generators\\
\hline
1&$(\frac{2}{p})=1$, $(\frac{b}{p})_4=1$&$\ge 5$&$\al$&$P_1P_2$& $P_1=(p,\al^2+u\al+t)$, $t^2=b\,\md{p}$  \\
& & && &$P_2=(p,\al^2-u\al+t)$, $2t=u^2\,\md{p}$ \\
\hline
2&$(\frac{-1}{p})=1$,$(\frac{-b}{p})_4=1$&$\ge5$&$\al$&$P_1P_2P_3P_4$&$P_1=(p,\al+t)$,
$P_2=(p,\al-t)$, $t^4=-b\,\md{p}$  \\
& & && &$P_3=(p,\al+ut)$, $P_4=(p,\al-ut)$ , $u^2=-1\,\md{p}$ \\
\hline
3&$(\frac{-1}{p})=-1$, $(\frac{-b}{p})_4=1$&$\ge
5$&$\al$&$P_1P_2P_3$&
$P_1=(p,\al+t)$, $P_2=(p,\al-t)$, $t^4=-b\,\md{p}$  \\
& & & &&$P_3=(p,\al^2+t^2)$ \\
\hline
4&$(\frac{-b}{p})=1$, $(\frac{-b}{p})_4\neq
1$&$\ge5$&$\al$&$P_1P_2$&$P_1=(p,\al^2+t)$,
$P_2=(p,\al^2-t)$, $t^2=-b\,\md{p}$  \\
\hline
5& $(\frac{4b}{p})_4=1$ ,$(\frac{-b}{p})_4\neq
1$&$\ge5$&$\al$&$P_1P_2$&$P_1=(p,\al^2+u\al+t)$,
$P_2=(p,\al^2-u\al+t)$  \\
& & & && $u^4=4b\,\md{p}$, $t^2=b\,\md{p}$ \\
\hline
6&$(\frac{-b}{p})_4\neq1$ &$\ge 5$&$\al$&$P$& $P=(p)$ \\
 &$(\frac{4b}{p})_4\neq 1$, $(\frac{-b}{p})=-1$& & & &  \\
\hline
7&$b=1\,\md{3}$&3&$\al$&$P_1P_2$&$P_1=(3,\al^2+\al-1)$  \\
& & & &&$P_1=(3,\al^2-\al-1)$\\
\hline
8&$b=-1\,\md{3}$&3&$\al$&$P_1^2P_2P_3$&$P_1=(3,\al^2+1)$  \\
& & & &&$P_2=(3,\al-1)$, $P_3=(3,\al+1)$\\
 \hline
9&$b=1\md{4}$, $a=0\md{4}$&2 &$\al$&$P^4$ & $P=(2,\al-1)$\\
 \hline
10&$b=3\md{4}$, $a=2\md{4}$&2 &$\al$&$P^4$ & $P=(2,\al-1)$\\
 \hline
 11&$b=1\md{4}$, $a=2\md{4}$&2 &&$P_1^3P_2$ & go to $table 8$\\
 \hline
 12&$b=7\md{8}$, $a=4\md{8}$&2 &&$P^2$ & $P=(2,\al-1)$\\
 \hline
 13&$b=7\md{8}$, $a=0\md{8}$&2 &$\frac{\t^3+4\t^2+6\t}{4}$&$P_1^2P_2$&$P_1=(2,\phi+1)$ \\
  &$1+b+a=8\md{16}$&&$\t=\al-1$&&$P_2=(2,\phi^2+\phi+1)$\\
  \hline
 14&$b=7\md{8}$, $a=0\md{8}$&2 &&$P_1^2P_2P_3$& go to $Table C14$ \\
  &$1+b+a=0\md{16}$&&&&\\
 \hline
15&$b=3\md{8}$, $a=4\md{8}$&2& & go to $Table C14$&\\
 \hline
 \end{tabular}$$}
{\small
$$Table D:\, v_p(ab)=0$$
For $p\ge5$, let $s\in \z$ such that
$3as+4b=0\md{p^{v_p(\triangle)}}$,
  and let $\theta=\al-s$.
$$\begin{tabular}{|c|c|c|c|c|c|
}\hline
Case&Conditions& $p$ &$\phi$&$p\z_K$& Generators\\
\hline
D1&&2& $\al$&$P$& $P=(2)$\\
\hline
D2&$b=-1\,\md{3}$&3& $\al$&$P$& $P=(3)$\\
\hline
D3&$a=b=1\,\md{3}$&3& $\al$&$P_1P_2$& $P_1=(3,\al-1)$, $P_2=(3,\al^3+\al^2+\al-1)$\\
\hline
D4&$a=b=1\,\md{3}$&3& $\al$&$P_1P_2$& $P_1=(3,\al+1)$, $P_2=(3,\al^3-\al^2+\al+1)$\\
\hline
D5&$ v_p(\triangle)=0$ &$p\ge 5$&$\al$ &$p$-analogous to$\bar P(X)$&\\
 \hline
D6&$ v_p(\triangle)=1$ &$p\ge 5$&$\theta$&   $P_1P_2P_3^2$  &$P_3=(p,\theta)$, $P_1=(p,\theta-v_1)$, $P_2=(p,\theta-v_2)$ \\
& $(\frac{-2}{p})=1$&   & & &$u^2=-2\,\md{p}, v_1=-s(2+u), v_2=-s(2-u)\,\md{p}$\\
 \hline
D7&$ v_p(\triangle)=1$ &$p\ge 5$&$\theta$&$P_1P_2^2$ &$P_2=(p,\theta)$,  \\
 &$(\frac{-2}{p})=-1$& & &&$P_1=(p,\theta^2+4s\theta+6s^2)$\\
 \hline
 D8&$ v_p(\triangle)\ge 2$ &$p\ge 5$& && go to $Table D8$\\
 \hline
\end{tabular}$$}
{\small  $$Table B6: B=a^4-a^2+b,\, A=-4a^3+a$$
$$\begin{tabular}{|c|c|c|c|c|c|c|}\hline
Case&Conditions& $3\z_K$&$\beta_3$& $\beta_2$ & $\beta_1$\\
\hline
B6.1& $ v_3(B)= 2$, $v_3(A)=1$&$P_1P_2^2P_3$&
$\frac{\theta^3-4a\t^2}{3}$ &
$\frac{\theta^3-4a\t^2+6sa^2\t+A}{3}+\t^2$   &$\theta-4a$\\
  \hline
B6.2& $v_3(B)= 3+k$, $v_3(A)=1$&$P_1P_2^2P_3$&
 $\frac{\theta^3-4a\t^2}{3}+3$ &
 $\frac{\theta^3-4a\t^2+6a^2\t+A}{3}+\t$   &$\theta-4a$\\
  \hline
 \end{tabular}$$}
{\small  $$Table  B11.4: B=s^4+as+b,\, A=4s^3+a,\,
as=-4b_3\md{3^{v_3(\triangle)+1}}$$
$$\begin{tabular}{|c|c|c|c|c|c|c|}\hline
Case&Conditions& $3\z_K$&$\beta_4$&$\beta_3$& $\beta_2$ & $\beta_1$\\
\hline
 1&$ v_3(b)\ge 2$, $a^2\neq 7\md{9}$&$P_1P_2^2P_3$&&
$\frac{\theta^3-4a\t^2}{3}$ &
$\frac{\theta^3-4a\t^2+6sa^2\t+A}{3}+\t^2$   &$\theta-4a$\\
 &$v_3(a^4-a^2+b)=2$&  &&  &  & \\
  \hline
2&$ v_3(b)\ge 2$, $a^2\neq 7\md{9}$, &$P_1P_2^2P_3$&
&$\frac{\theta^3-4a\t^2}{3}+3$ &
 $\frac{\theta^3-4a\t^2+6a^2\t+A}{3}+\t$   &$\theta-4a$\\
 &$v_3(a^4-a^2+b)\ge 3$&  &                          & &  & \\
 \hline
3&$ v_3(b)\ge 2$, $a^2= 7\md{9}$, &$P_1P_2^3$&  && $\frac{\theta^3+4s\t^2}{3}+\t$   &$\theta+4s$\\
 &$v_3(a^4-a^2+b)=2$&  &&  &  & \\
 \hline
4&$ v_3(\triangle)=2r+1$, $r\ge 4$ &$P_1P_2P_3^2$&
&$\frac{\theta^3+4s\t^2+6s^2\t+A}{3^{r-1}}+\frac{\theta^2+4s\t}{3}$
& $\frac{\theta^3+4s\t^2}{3}+\t$ &$\theta+4s$\\
 \hline
5&$ v_3(\triangle)=7$,  &$P_1P_2P_3^2$&&
$\frac{\theta^3+4s\t^2+6s^2\t+A}{3^{r-1}}+\frac{\theta^2+4s\t}{3}$
& $\frac{\theta^3+4s\t^2}{3}+\t^2$   &$\theta+4s$\\
  \hline
6& $ v_3(\triangle)=2r$, $r\ge 5$ &$P_1P_2P_3$&&
 $\frac{\theta^2+4s\t}{3}+\t$ & $\frac{\theta^2+4s\t+6s^2}{3}+\t$
&$\theta+4s$\\
  &         $(\frac{-2B_3}{3})=-1$                      &  &&  &  & \\
 \hline
 7& $ v_3(\triangle)=2r\ge 8$&$P_1P_2P_3P_4$&$\frac{\t^3+4s\t^2+6s^2(\theta+3^{r-2}t)}{3^{r-1}}$&
$\frac{\t^3+4s\t^2+6s^2(\theta-3^{r-2}t)}{3^{r-1}}$&$\frac{\theta^2+4s\t+6s^2}{3}+\t$
&$\theta+4s$\\
   &         $(\frac{-2B_3}{3})=1$ && $2s^2t^2+B_3=3\md{9}$ &  &  & \\
  \hline
\end{tabular}$$}
 If ($b=7\md{8}$, $a=0\md{8}$ and $1+b+a=0\md{16}$) or $b=3\md{8}$ and  $a=4\md{8}$,
 then let $s\in\z$ such that $P(X)$ is $(X+s)$-regular, and let $A=4s+a$, $B=s^4+as+b$ and  $\theta=\al-s$.
  {\small
$$Table C14: v_2(B)\ge 3, v_2(A)\ge 2$$
$$\begin{tabular}{|c|c|c|c|c|c|c|c|}\hline
Case&Conditions& $2\z_K$&$\beta_3$& $\beta_2$& $\beta_1$   \\
\hline
1& $ v_2(B)\ge  2$ &$P_1^3P_2$&&$\frac{\theta^3}{2}+2$ & $\frac{\t^3+A}{2}+\t$  \\
 &$v_2(A)=1$&  & & &  \\
 \hline
2&$ v_2(B)= 2r-1$ &$P_1^2P_2$& &$\frac{\t^2}{2}+2$& $\frac{\theta^2+6}{2}+\t$ \\
 &$v_2(A)=r\ge 2$&  &  & &\\
   \hline
 3&$ v_2(B)> 2r$ &$P_1^2P_2P_3$&
 $\frac{\t^3+4s\t^2+6s^2\t}{2^r}$ & $\frac{(\t^2+4s\t+6s^2)(\t+2^{r-1}t)}{2^r}+2$
  &$\frac{\theta^2+6}{2}+\t$\\
 &$v_2(A)=r\ge 2$&  &  & $A_2t=1\md{4}$&\\
 \hline
 4&$ v_2(B)= 2r$ &$P_1^2P_2P_3$&
 $\frac{\t^3+4s\t^2+6s^2\t}{2^r}$ & $\frac{(\t^2+4s\t+6s^2)(\t+2^{r-1}t)}{2^r}$
  &$\frac{\theta^2+6}{2}+\t$\\
 &$v_2(A)=r\ge 2$&  &  & $A_2t=3\md{4}$&\\
 \hline
5&$ v_2(B)= 2r$ &$P_1^2P_2^2$&  &$\frac{\theta^2}{2}+\frac{\t^3+4s\t^2+6s^2\t+A}{2^r}$ & $\frac{\theta^2+6}{2}+\t$  \\
 &$v_2(A)\ge r+1$&  &  &  &\\
  \hline
\end{tabular}$$}
{\small   $$Table D8 : p\ge 5, \, v_p(ab)=0 {\mbox and }
v_p(\triangle)\ge 2$$ Let  $s\in\z$ such that
$3as+4b=0\md{p^{v_p(\triangle)+1}}$, $\theta=\al-s$, $B=s^4+as+b$
and $A=4s^3+a$.
{\small$$\begin{tabular}{|c|c|c|c|c|c|c|c|}\hline
Conditions& $p\z_K$&$\beta_4$&$\beta_3$& $\beta_2$& $\beta_1$\\
\hline
$ v_p(\triangle)=2r+1$ & $P_1P_2P_3^2$& &  $\frac{\theta^3+4s\t^2+6s^2\t}{p^{r}}$ & $\theta-v_1$& $\theta-v_2$\\
$(\frac{-2}{p})=1$ &  & & $B_p=6s^2u\md{p}$ &$u^2+2=p\md{p^2}$& $u^2+2=p\md{p^2}$\\
& &   && $v_1=-s(2+u)\md{p^2}$&$v_2=-s(2-u)\md{p^2}$ \\
 \hline
$ v_p(\triangle)=2r+1$ &$P_1P_2^2$&  & & $\frac{\theta^3+4s\t^2+6s^2\t}{p^{r}}$& $\theta^2+4s\theta+6s^2+p$ \\
$(\frac{-2}{p})=-1$& & &  & $B_p=6s^2u\,\md{p}$&\\
 \hline
$ v_p(\triangle)=2r$ &$P_1P_2P_3$& & $\frac{\theta^3+4s\t^2+6s^2\t}{p^{r-1}}$  &$\theta-v_1$&$\theta-v_2$\\
 $(\frac{-2}{p})=1$& && $v_p(6s^2+B_p)=1$ & $u^2+2=p\md{p^2}$& $u^2+2=p\md{p^2}$\\
 $(\frac{-6B_p}{p})=-1$ & &  & & $v_1=-s(2+u)\md{p^2}$&$v_2=-s(2-u)\md{p^2}$ \\
  \hline
 $ v_p(\triangle)=2$ &$P_1P_2P_3P_4$& $\frac{(\t^2+4s\t+6s^2)(\theta+pt)}{p}+\t^2$ & $\frac{(\t^2+4s\t+6s^2)(\theta-pt)}{p}+\t^2$ &$\theta-v_1$&$\theta-v_2$\\
$(\frac{-2}{p})=1$& &$v_p(6s^2t^2+B_p)\ge 2$&  & $u^2+2=p\md{p^2}$& $u^2+2=p\md{p^2}$\\
$(\frac{-6B_p}{p})=1$& &  & & $v_1=-s(2+u)\md{p^2}$&$v_2=-s(2-u)\md{p^2}$ \\
  \hline
  $ v_p(\triangle)=2r\ge 4$ &$P_1P_2P_3P_4$& $\frac{(\t^2+4s\t+6s^2)(\theta+p^{r}t)}{p^r}+\t$ & $\frac{(\t^2+4s\t+6s^2)(\theta-p^{r}t)}{p^r}+\t$  &$\theta-v_1$&$\theta-v_2$\\
$(\frac{-2}{p})=1$& &$v_p(6s^2t^2+B_p)=1$& & $u^2+2=p\md{p^2}$& $v_p(u^2+2)=1$\\
$(\frac{-6B_p}{p})=1$& &  & & $v_1=-s(2+u)\md{p^2}$&$v_p(v_2+s(2-u))\ge2$ \\
  \hline
 $ v_p(\triangle)=2r$ &$P_1P_2$&  & & $\frac{\theta^3+4s\t^2+6s^2\t}{p^{r-1}}$ & $\theta^2+4s\theta+6s^2+p$ \\
 $(\frac{-2}{p})=-1$& & &  & &\\
 $(\frac{-6B_p}{p})=-1$& &  & &&\\
  \hline
 $ v_p(\triangle)=2$ &$P_1P_2P_3$& &$\frac{(\t^2+4s\t+6s^2)(\theta+pt)}{p}+\t^2$ & $\frac{(\t^2+4s\t+6s^2)(\theta-pt)}{p}+\t^2$& $\theta^2+4s\theta+6s^2+p$ \\
 $(\frac{-2}{p})=-1$& & &$v_p(6s^2t^2+B_p)\ge2$&  & $v_p(6s^2u+B_p)=1$\\
 $(\frac{-6B_p}{p})=1$& &  & &&\\
  \hline
 $ v_p(\triangle)=2r\ge 4$ &$P_1P_2P_3$& &$\frac{(\t^2+4s\t+6s^2)(\theta+p^{r}t)}{p^r}+\t$ & $\frac{(\t^2+4s\t+6s^2)(\theta-p^{r}t)}{p^r}+\t$  & $\theta^2+4s\theta+6s^2+p$ \\
 $(\frac{-2}{p})=-1$& & &$v_p(6s^2t^2+B_p)=1$  & &\\
  $(\frac{-6B_p}{p})=1$& &  & &&\\\hline
\end{tabular}$$}}
\end{teor}
{\bf Proof of Theorem}. All cases, except ($p=2$, $v_p(b)=2$ and
$v_p(a)\ge 2$), correspond to a situation where $P(X)$ is
$p$-regular. The case : $p=2$, $v_p(b)=2$ and $v_p(a)\ge 2$ is
handled in $Table A8$ by using technics of  Newton polygons of
second order. \\
Denote $C(X)$ the minimal polynomial of a possible
$\phi$ such that  $v_p([\z_K:\z[\phi]])=0$. \\
 {\bf $v_p(a)\ge 1$ and $v_p(b)\ge 1$.}
  \begin{enumerate}
    \item $v_p(a)\ge 3$ and $v_p(b)=3$. Then $\bar P(X)=X^4\,\md{p}$
    and the $X$-Newton polygon of $P(X)$ is one side of slope
    $3/4$. Thus, $p\z_K=P^4$. Since $v_P(\al)=4\times 3/4=3$, then
    $v_P(\frac{\alpha^3}{p^2})=1$ and
    $P=(p,\frac{\alpha^3}{p^2})$.
\item  Let $\phi=\frac{\alpha^3}{p^h}$. Then
$C(X)=X^4+3\frac{a}{p^h}X^3+3\frac{a^2}{p^{2h}}X^2+\frac{a^3}{p^{3h}}X+\frac{b^3}{p^{4h}}$.
It follows that if $v_p(a)=1$ and $v_p(b)\ge 2$, then
$\phi=\frac{\alpha^3}{p}$ is a $p$-generator of $\z_K$, $\bar
C(X)=X(X+a_p)^3\,\md{p}$ and $p\z_K=(p,\al)(p,\al+a_p)^3$. If
$v_p(a)=2$ and $v_p(b)\ge 3$, then $\phi=\frac{\alpha^3}{p^2}$ is
a $p$-generator of $\z_K$, $\bar C(X)=X(X+a_p)^3\,\md{p}$ and
$p\z_K=(p,\phi)(p,\phi+a_p)^3$.
 \item If $v_p(b)=1$ and $v_p(a)\ge 1$, then $v_p(ind(P))=0$,
 then $\al$ is a $p$-generator of $\z_K$ and $p\z_K=(p,\al)^4$.
 \item
 If $v_p(b)=2$ and $v_p(a)\ge 2$, then $N_X(P)=S$ is one side  such
 that $P_S(y)=Y^2+b_p$. It follows that:
\begin{enumerate}
\item
 If $(\frac{-b_p}{p})=-1$, then $p\z_K=P^2$ such that
 $v_P(\al)=1$. Thus, $P=(p,\frac{\al^3}{p})$.
 \item
If $(\frac{-b_p}{p})=1$ and  $v_p(a)=2$, then $p\z_K=P_1^2P_2^2$.
For $\phi=\frac{\al^2}{p}$, we have
$C(X)=X^4-Xpa_p^2+2b_pX^2+b_p^2$ and  $\bar
C(X)=(X^2+b_p)^2\,\md{p}$. Let $t\in\z$ such that $v_p(t^2+b_p)=1$.
Let $\ph(X)=X+t$ and
$\ph(X)^4-4t\ph(X)^3+(6t^2+2B)\ph(X)^2+(-pa_p^2-4tb_p-4t^3)\ph(X)+(b_p^2+tpa_p^2+2t^2b_p+t^4)$
be the $\ph(X)$-adic development of $C(X)$. Since
$v_p(b_p^2+tpa_p^2+2t^2b_p+t^4)=1$, $\frac{\al^2}{p}$ is a
$p$-generator of $\z_K$. As $\bar C(X)=(X^2+b_p)^2\, \md{p}$, we have
$P_1=(p,t+\frac{\al^2}{p})$ and $P_2=(p,-t+\frac{\al^2}{p})$
$v_p(t^2+b_p)=1$. \item
 If $(\frac{-b_p}{p})=1$ and  $v_p(a)=\ge 3$, then $p\z_K=P_1^2P_2^2$.
 Let $\phi=\frac{\alpha^3+\al^2}{p}$.
 Since $\alpha^3+\al^2\not\in \q$, $K=\q[\alpha^3+\al^2]$  and
$\bar C(X)= X^4+2b_pX^2-4pb_p^2X+b_p^2\md{p^2}$ and  $\bar
C(X)=(X^2+b_p)^2\,\md{p}$. Let $\phi(X)=X+t$ such that $t^2+b_p=0\,
\md{p}$ and $C(X)=\phi^4(X) -4t\phi^3(X)+(6t^2+2b_p)\phi^2(X)+(
-4pb_p^2-4tb_p-4t^3)\phi(X)+(b_p^2+4tpb_p^2+2t^2b_p+t^4)$ be the
$\ph(X)$-adic development of $C(X)$. Since
$v_p(b_p^2+4tpb_p^2+2t^2b_p+t^4)=1$, $\frac{\al^3+\al^2}{p}$ is a
$p$-generator of $\z_K$, $P_1=(p,t+\frac{\al^3+\al^2}{p})$ and
$P_2=(p,-t+\frac{\al^3+\al^2}{p})$.
\end{enumerate}
 \item
$v_2(b)= 2$ and $v_2(a)\ge 2$. In that case, using technics of
Newton polygons of second order, for every subcase, a monic
polynomial $\phi_2\in\z[X]$ such that $t=(X, 1/2, \phi_2)$ is
$P(X)$-complet is given (see \cite [Def 3.9, p 38]{GNM}). Let
$N_2$ be the $\phi_2$-Newton polygon of second order of $P(X)$.
From \cite [Cor 3.8, p 38]{GNM}, if $N_2$ is one side, then
$2\z_K=P^e$, where $e=2e_2$ and $e_2$ is the ramification index of
$N_2$. If $N_2$ is two sides, then for every side, $e_2=1$ and
then $2\z_K=P_1^2P_2^2$.
\begin{enumerate}
\item
 If $v_2(a)= 2$, then for $\phi_2=X^2+2$ and $N_2$ is
one side. Thus, $2\z_K=P^4$,  $v_P({\al^2+2})=5$ and
$2\z_K=(2,\frac{\al^2+2}{2})$.
 \item
 If $v_2(a)= 3$ and $b=4\md{16}$ (resp. $v_2(a)= 3$ and  $b=12\md{16}$), then for
$\phi_2=X^2+2X+2$ (resp. $\phi_2=X^2+2$) and  $N_2$ is one side.
Then $2\z_K=P^4$, $v_P({\al^2+2\al+2})=7$ (resp.
$v_P({\al^2+2})=7$) and $2\z_K=(2,\frac{\al^3+2\al^2+2\al}{4})$
(resp. $2\z_K=(2,\frac{\al^3+2\al}{4})$).
 \item
  If $v_2(a)= 4$ and $b=4\md{32}$, then for
$\phi_2=X^2+2X+2$,  $N_2$ is one side, $2\z_K=P^4$ and
$v_P({\al^2+2\al+2})=9$. Thus, $2\z_K=(2,\frac{\al^2+2\al+2}{4})$.
\item
  If $v_2(a)= 4$ and $b=20\md{64}$, then for
$\phi_2=X^2+2X-2$,  $N_2$ is one side, $2\z_K=P^2$ and
$v_P({\al^2+2\al-2})=5$. Thus,
$2\z_K=(2,\frac{\al^2+2\al-2\al}{4})$.
 \item
  If $v_2(a)\ge 4$ and $b=12\md{32}$, then for
$\phi_2=X^2+2$,  $N_2$ is one side, $2\z_K=P^2$ and
$v_P({\al^2+2})=4$. Thus, $2\z_K=(2,\frac{\al^3+2\al}{4})$.
\item
  If $v_2(a)\ge 5$ and $b=20\md{32}$, then for
$\phi_2=X^2+2X-2$,  $N_2$ is one side, $2\z_K=P^4$ and
$v_P({\al^2+2\al-2})=9$. Thus,
$2\z_K=(2,\frac{\al^2+2\al-2\al}{4})$.
\item
  If $v_2(a)\ge5$ and $b=36\md{64}$, then for
$\phi_2=X^2+2X+2$,  $N_2$ is one side, $2\z_K=P^2$ and
$v_P({\al^2+2\al+2})=5$. Thus,
$2\z_K=(2,\frac{\al^2+2\al-2\al}{4})$.
\end{enumerate}
\item $b=4+64B$ and $a=32+64A$. For
$\phi=\frac{\alpha^3-2\al-4}{8}$, we have $\bar
C(X)=X^4+2X^3+3X^2+2X+2\md{4}$. Thus, $v_2(C(1))=v_2(C(0))=1$,
$\phi$ is a $2$-generator of $\z_K$ and
$2\z_K=(2,\phi)^2(2,\phi+1)^2$.
 \item $b=4+64B$ and $a=64A$. For
$\phi=\frac{\alpha^3+4\al^2-2\al+4}{8}$, we have $\bar
C(X)=X^4+2X^3+3X^2+2X+2\md{4}$. It follows that
$v_2(C(1))=v_2(C(0))=1$, $\phi$ is a $2$-generator of $\z_K$  and
$2\z_K=(2,\phi)^2(2,\phi+1)^2$.
 \item $b=28+32B$ and  $a=16+32A$.
For $\phi=\frac{\alpha^2+2}{4}$, we have $C(X)=
X^4-2X^3+(5+4B)X^2+(-16A^2-8-16A-4B)X+6+8A+8A^2+8B+4B^2=(X+1)^2X^2\,\md{2}$.
Since $v_2(C(0))=v_2(C(1))=1$, we  have $\phi$ is a $2$-generator
of $\z_K$ and $2\z_K=(2,\phi+1)^2(2,\phi)^2$.
 \item
 $b=52+64B$ and $a=16+64A$. For $\phi=\frac{\alpha^3-2\al^2-2\al}{8}$,  $C(X)=
X^4+2X^3+3X^2+2X+2\md{4}$. Thus, $v_2(C(0))=v_2(C(1))=1$, $\phi$ is
a $2$-generator of $\z_K$ and $2\z_K=(2,\phi+1)^2(2,\phi)^2$.
\item
 $b=52+64B$ and $a=48+64A$. For $\phi=\frac{\alpha^3-2\al^2+6\al}{8}$,  $C(X)=
X^4+2X^3+3X^2+2X+2\md{4}$. Thus, $v_2(C(0))=v_2(C(1))=1$, $\phi$ is
a $2$-generator of $\z_K$ and $2\z_K=(2,\phi+1)^2(2,\phi)^2$.
 \item  $b=28+32B$ and  $a=32A$.
For $\phi=\frac{\alpha^2+12\al+2}{4}$,  $C(X)=
X^4+2X^3+X^2+2\md{4}$. Thus, $v_2(C(0))=v_2(C(1))=1$, $\phi$ is a
$2$-generator of $\z_K$ and $2\z_K=(2,\phi+1)^2(2,\phi)^2$.
\end{enumerate}
 {\bf $v_p(a)=0$ and   $v_p(b)\ge 1$.}\\
  If $p\neq 3$, then $v_p(ind(P))=0$, $\al$ is a
 $p$-genrator of $\z_K$  and $p\z_K$ is $p$-analogous to $\bar
 P(X)$.\\
     For $p=3$, let
$F(X)=P(X-a)=X^4-4aX^3+6a^2X^2+AX+B$ and $\t=\al+a$
($A=-a(4a^2-1)$ and $B=(a^4-a^2+b)$).
\begin{enumerate}
\item If $v_3(B)=1$,  then $v_3(ind(P))=0$ and
$3\z_K=(3,\al+a)^3(3,\al)$.
\item If $v_3(B)=2$ and $v_3(A)\ge 2$,
then $3\z_K=P_1P_2^3$, $v_{P_1}(\t-4a)=1$, $v_{P_1}(\t)=0$,
$v_{P_2}(\t-4a)=0$ and $v_{P_2}(\t)=2$. Thus, $P_1=(3,\t-4a)$,
$P_2=(3,\frac{\t^3-4a\t^2}{3})$.
 \item If $v_3(A)=1$ and
$v_3(B)\ge 2$, then $3\z_K=P_1P_2^2P_3$.\\
Since $\t^4-4a\t^3=-6a^2\t^2+A\t+B$ and
$v_{P_1}(\t)=0$, we have $v_{P_1}(\t^3-4a\t^2)=1$ and $P_1=(3,\t-4a)$.\\
If $v_3(B)= 2$, then let $\beta_3=\frac{\theta^3-4a\t^2}{3}$ and
$\beta_2=\frac{\theta^3-4a\t^2+6sa^2\t+A}{3}+\t^2$. Since
 $v_{P_3}(\t)=v_{P_2}(\t)=1$ and $v_{P_1}(\t)=0$, we have
$v_{P_3}(\beta_3)=1$ and $v_{P_i}(\beta_3)=0$ for $i\neq 3$. We
have also, $v_{P_2}(\beta_2)=1$ and $v_{P_i}(\beta_3)=0$ for
$i\neq 2$.\\
The case $v_3(B)\ge 3$ is similar to the previous case.
 \item
 If $v_3(A)\ge 2$ and $v_3(B)\ge 3$ ($v_3(\triangle)\ge 6$), then  let $s\in \z$ such that
$as=-4b_3\md{3^{v_3(\triangle)+1}}$. Let
$F(X)=P(X+s)=X^4+4sX^3+6s^2X^2+AX+B$ and  $\theta=\al-s$. Since
$A=4s^3+a$ and $B=s^4+as+b$, $v_3(A)=v_3(B)=v_3(\triangle)-3\ge
3$. It follows that:\\
    If $v_3(\triangle)=6$, then for
       $\phi=\frac{\theta^3+4s\theta^2+6s^2\theta+A}{9}$,
  we have $\bar C(X)=X(X^3+2s^2B_3X-4sB_3^2)\,\md{3}$. So, if  $B_3=1\,\md{3}$,
  then $\bar C(X)=X(X^3-X-s)\,\md{3}$ and $3\z_K=(3,\phi)(3,\phi^3-\phi-s)$.
  If $B_3=-1\,\md{3}$ and $s=1\,\md{3}$, then $\bar C(X)=X(X^2-X-1)(X+1)\,\md{3}$.
  If $B_3=-1\,\md{3}$ and $s=-1\,\md{3}$, then $\bar C(X)=X(X^2+X-1)(X-1)\,\md{3}$
\item
 If $v_3(\triangle)\ge 7$,
 then $N_X(F)=S_3+S_2+S_1$, $F_{S_1}(Y)=Y+4s$ and  $F_{S_2}(Y)=sY-1$. Since
$\t^3(\t^2+4s)=-(6s^2\t^2+A\t+B)$, then  $P_1=(3,\theta+4s)$.
\begin{enumerate}
\item
   If $v_3(\triangle)=2r+1$ ($r\ge 3$), then   $v_3(A)=v_3(B)=2r-2$, $F_{S_3}(Y)=2s^2Y+B_3$,
    $3\z_K=P_1P_2P_3^2$,  $v_{P_3}(\theta)=2r-3$, $v_{P_2}(\theta)=1$ and
    $v_{P_1}(\theta)=0$. For every $i\neq j$, $v_{P_i}(\beta_i)=1$ and
    $v_{P_i}(\beta_j)=0$.
 \item If $ v_3(\triangle)=2r$,
then $F_{S_3}(Y)=6s^2Y^2+B_3$. Thus, if $(\frac{-2B_3}{3})=-1$,
then $3\z_K=P_1P_2P_3$. If $(\frac{-2B_3}{3})=-1$, then $3\z_K=P_1P_2P_3P_4$.\\
 If $(\frac{-2B_3}{3})=-1$, then let  $\beta_3=\frac{\theta^2+4s\t}{3}+\t+3$  and
 $\beta_2=\frac{\theta^2+4s\t+6s^2}{3}+\t$.\\
If  $(\frac{-2B_3}{3})=1$, then $F_{S_3}(Y)=6s^2(Y+t)(Y-t),$ where
$6s^2t^2+B_3=0\,\md{3}$. Thus, $3\z_K=P_1P_2P_3P_4$, where $P_3$ and
$P_4$ are respectively attached to $Y+t$ and $Y-t$.  Let
$\phi=X+3^{r-2}t$ and  consider the $\phi$-adic development of
$F(X)$, where $2s^2t^2+B_3=3\md{9}$. We have
$v_{P_3}(\phi(\theta))=r-1$ and $v_{P_4}(\phi(\theta))=r-2$.
Therefore,
$\beta_4=\frac{\t^3+4s\t^2+6s^2(\theta+3^{r-2}t)}{3^{r-1}}$,
$\beta_3=\frac{\t^3+4s\t^2+6s^2(\theta-3^{r-2}t)}{3^{r-1}}$ and
$\beta_2=\frac{\theta^2+4s\t+6s^2}{3}+\t$, where
$2s^2t^2+B_3=3\md{9}$ satisfy:
$v_{P_3}(\t^3+4s\t^2)=v_{P_4}(\t^3+4s\t^2)=2(r-1)$,
$v_{P_3}(\theta+3^{r-2}t)=r-1$, $v_{P_4}(\theta+3^{r-2}t)=r-2$,
$v_{P_3}(\theta-3^{r-2}t)=r-2$ and $v_{P_4}(\theta-3^{r-2}t)=r-1$.
So, $v_{P_i}(\beta_i)=1$ and  $v_{P_j}(\beta_i)=0$ for every
$i\neq j$.
\end{enumerate}
{\bf $v_p(a)\ge 1$ and $v_p(b)=0$.}\\
 If $p\neq 2$, then
$v_p(ind(P))=0$. Thus, $\alpha$ ia $p$-generator of $\z_K$,  and
 then $p\z_K$ is $p$-analogous to $\bar P(X)=X^4+aX+b\,\md{p}$.
For $p\ge 5$, if $(\frac{-b}{p})_4=1$, then $\bar
P(X)=(X-u)(X+u)(X^2+u^2)$, where $u^4=-b\,\md{p}$.\\
$(\frac{-b}{p})_4\neq 1$, then $\bar P(X)$ is irreducible in
$\F_p[X]$ or  $\bar P(X)=(X^2+rX+s)(X^2-rX+t)\,\md{p}$. In the last
case, ($r=0\,\md{p}$, $s=-t\,\md{p}$ and $s^2=-b\,\md{p}$) or ($r\neq
0\,\md{p}$, $s=t\,\md{p}$, $2s=r^2\,\md{p}$ and $s^2=b\,\md{p}$); i.e., $\bar
P(X)=(X^2+s)(X^2-s)\,\md{p}$, where $s^2=b\,\md{p}$ or $\bar
P(X)=(X^2+uX+s)(X^2-uX+s)\,\md{p}$, where $u^4=4b\,\md{p}$ and
$2s=u^2\,\md{p}$.\\
 For $p=2$, let $F(X)=P(X+1)=X^4+4x^3+6X^2+AX+B$ and
$\t=\al-1$.
\begin{enumerate}
\item
 If $v_2(B)=3$, $A=4\md{8}$, then for
$\phi=\frac{\t^3+4\t^2+6\t}{4}$, we have
$C(X)=X^4+(3+2K)X^3+(2+2L)X^2+(3+2K)X+1+2L\md{4}$, where $A=4+8K$
and $B=8+16L$. Since $C(1)=2\md{4}$, $\phi$ is a $2$-generator of
$\z_K$,  $2\z_K=P_1^2P_2$, where $P_1=(2,\phi+1)$  and
$P_2=(2,\phi^2+\phi+1)$.
 \item
 $b=3\md{8}$ and $a=4\md{8}$ or ($b=7\md{8}$, $a=0\md{8}$ and
 $1+b+a=0\md{16}$), then let $s\in\z$ such that $P(X)$ is $X+s$-regular.
  Let $F(X)=P(X-s)=X^4+4sX^3+6s^2X^2+AX+B$ and  $\theta=\al-s$.
 It follows that:
\begin{enumerate}
\item $v_2(A)=1$ and $v_2(B)\ge 2$. Then $2\z_K=P_1^3P_2$.
  \item   If $v_2(A)=r$ and $v_2(B)=2r-1$, then $N_X(F)=S_1+S_2$
   such that $F_{S_1}(Y)=Y+1$ and $F_{S_2}(Y)=Y^2+Y+1$. Thus,
   $2\z_K=P_1^2P_2$, where  $v_{P_1}(\theta)=1$ and $v_{P_2}(\theta)=r-1$.
   Hence, for every $i\neq j$, $v_{P_i}(\beta_i)=1$  and $v_{P_i}(\beta_j)=0$.
\item
If $v_2(A)=r$ and $v_2(B)\ge  2r$, then $N_X(F)=S_1+S_2+S_3$
   such that $F_{S_1}(Y)=Y+1$ and $F_{S_3}(Y)=F_{S_2}(Y)=Y+1$.
   Thus,  $2\z_K=P_1^2P_2P_3$. Let
  $\beta_3=\frac{\t^3+4s\t^2+6s^2\t}{2^r}$,
 $\beta_2=\frac{(\t^2+4s\t+6s^2)(\t+2^{r-1}t)}{2^r}$
   and $\beta_1=\frac{\theta^2+6}{2}+\t$, where $t\in\z$
   is choused such that $v_{P_2}(\t+2^{r-1}t)=r$:  Let $\phi(X)=X+2^{r-1}t$ and
   consider the $\phi(X)$-adic of $F(X)$. If $v_2(B)= 2r$, then
   $A_2t=3\md{4}$ and f $v_2(B)> 2r$, then $A_2t=1\md{4}$.
\item
If $v_2(A)= r+1+k$ and $v_2(B)=2r$, then $N_X(F)=S_1+S_2$
   such that $F_{S_1}(Y)=Y+1$ and $F_{S_3}(Y)=Y+1$.
   Thus, $2\z_K=P_1^2P_2^2$, where  $v_{P_1}(\theta)=1$ and $v_{P_2}(\theta)=2r-1$.
\end{enumerate}
\end{enumerate}
{\bf $v_p(ab)=0$.}\\
If  $p\in \{2,3\}$, then $v_p(ind(P))=0$ and then $p\z_K$ is
$p$-analogous to $\bar P(X)$. For $p\ge 5$, let $s\in \z$ such
that $3as+4b=0\md{p^{v_p(\triangle)+1}}$, $\theta=\al-s$ and
$F(X)=P(X+s)=X^4+4sX^3+6s^2X^2+AX+B$. Then $v_p(A) = v_p(B) =
v_p(\triangle)$.\\
  If $v_p(\triangle)=0$, then  $p\z_K$ is $p$-analogous to $\bar P(X)$.
  If $v_p(\triangle)=1$ and $(\frac{-2}{p})=1$, then  $p\z_K=P_1^2P_2P_3$,
  where $P_1=(p,\theta)$, $P_2=(p,\theta+u)$, $P_3=(p,\theta-u)$, $\theta =\al-s$
   and $3at+4b=0\md{p^2}$. \\
   If $v_p(\triangle)=1$ and $(\frac{-2}{p})=-1$, then $p\z_K=P_1^2P_2$, where
    $P_1=(p,\theta)$ and  $P_2=(p,\theta^2+4s\theta+6s^2)$.
  If $v_p(\triangle)\ge 2$, then let  $s\in\z$ such that
$3as+4b=0\md{p^{v_p(\triangle)+1}}$, $\theta=\al-s$ and
$F(X)=P(X+s)=X^4+4sX^3+6s^2X^2+AX+B$, where  $B=s^4+as+b$ and
$A=4s^3+a$. Then $v_p(A)=v_p(B)=v_p(\triangle)$ and
$N_{X}(F)=S_0+S_1$ with respective slopes $0$ and
$\frac{v_p(\triangle)}{2}$. It follows that:
  \begin{enumerate}
 \item
If $ v_p(\triangle)=2r+1$, then $N_X(F)=S_1+S_2$
   such that $F_{S_1}(Y)=Y^2+4sY+6s^2$ and $F_{S_2}(Y)=6s^2Y+B_p$.
   Thus, if $(\frac{-2}{p})=1$, then
$p\z_K=P_1P_2P_3^2$. If $(\frac{-2}{p})=-1$, then
$p\z_K=P_1P_2^2$.
 \item
If $ v_p(\triangle)=2r$, then $N_X(F)=S_1+S_2$
   such that $F_{S_1}(Y)=Y^2+4sY+6s^2$ and $F_{S_2}(Y)=6s^2Y^2+B_p$.
   Thus, $(\frac{-2}{p})=1$ and
$(\frac{-6B_p}{p})=-1$, then $p\z_K=P_1P_2P_3$.\\
If $(\frac{-2}{p})=-1$ and $(\frac{-6B_p}{p})=-1$, then
$p\z_K=P_1P_2$.\\
 If $(\frac{-2}{p})=-1$ and
$(\frac{-6B_p}{p})=1$, then $F_{S_2}(Y)=6s^2(Y-t)(Y+t)$, where $t$
is a root of $F_{S_2}(Y)$ in $\F_3$. Thus, $p\z_K=P_1P_2P_3$,
where $P_1$, $P_2$ and $P_3$ are respectively attached to $S_1$,
$Y+t$ and $Y-t$. Let
$\beta_2=\frac{(\t^2+4s\t+6s^2)(\theta+p^{r}t)}{p^r}+\t$ and
$\beta_3=\frac{(\t^2+4s\t+6s^2)(\theta-p^{r}t)}{p^r}+\t$, where
$t\in\z$ is choused such that $v_{P_2}(\theta+p^{r}t)=r+1$,
$v_{P_3}(\theta+p^{r}t)=r$, $v_{P_2}(\theta+-p^{r}t)=r$ and
$v_{P_3}(\theta-p^{r}t)=r+1$: If $r\ge 2$, then
$6s^2t^2+B_p=p\md{p^2}$. If $r=1$, then $6s^2t^2+B_p=0\md{p^2}$. In
that way, we have
$P_2=(p,\frac{(\t^2+4s\t+6s^2)(\theta+p^{r}t)}{p^r}+\t)$ and
$P_3=(p,\frac{(\t^2+4s\t+6s^2)(\theta-p^{r}t)}{p^r}+\t)$.\\
 If $(\frac{-2}{p})=1$ and $(\frac{-6B_p}{p})=1$, then
$p\z_K=P_1P_2P_3P_4$.
 \end{enumerate}
\end{enumerate}
\begin{exs}
Let $P(X)=X^4+aX+b\in\z[X]$ be an irreducible polynomial, $\al$ a complex root of $P(X)$ and $K=\q[\al]$.
\begin{enumerate}
\item
$a=2^{10}.5$ and $b=2^{9}.3.5$. Let $\t=\frac{\al}{4}$. Then $\t$ is integral with minimal polynomial $F(X)=X^4+80X+30$ and $K=\q[\t]$. Thus
$2\z_K=(2,\t)^4$, $3\z_K=(3,\t)(3,\t-1)^3$, $5\z_K=(5,\t)^4$ and for every $p\not\in\{2,3,5\}$, $p\z_K$ is $p$-analogous to $\bar F(X)$.
\item
$a=48$ and $b=188$. From $Table A8$, row $A8.9$, we have $2\z_K=(2,\frac{\al^2+2}{4})^2(2,\frac{\al^2+6}{4})^2$.
\item
$a=144$ and $b=36$. From $Table A8$, row $A8.4$, we have
$2\z_K=(2,\frac{\al^2+2\al+2}{4})^4$. Since
$(\frac{-b_3}{3})=-1$, from $Table A$, row $A4$ $3\z_K=P^2$,
where $P=(3,\frac{\al^3}{3})$. Since $\triangle=-2^{14}.3^6.971$,
for every $p\not\in\{2,3,5\}$, $p\z_K$ is $p$-analogous to $\bar
P(X)$.
\item
$a=28$, $b=189$ and $p=2$. Let $F(X)=P(X+1)=X^4+4X^3+6X^2+32X+218$. Since $B=218=2.109$, we have  $2\z_K=(2,\al-1)^4$.
\item
$a=22$ and $b=66$. Then $\triangle=2^4.3^5.11^3.13$. We have
$2\z_K=(2,\al)^4$, $11\z_K=(11,\al)^4$,
$13\z_K=(13,\al)(13,\al^3+9)$ and for every prime
$p\not\in\{2,3,11,13\}$, $p\z_K$ is $p$-analogous to $\bar P(X)$.
For $p=3$, since $a^2=7\md{9},\, b=3 \md{9}$ and
$b+a^4-a^2=9.25982$,  from $Table B$, row $B10$,
$3\z_K=P_1P_2^3$, where $P_1=(3,\al-5a)$ and
$P_2=(3,\frac{\t^3-4a\t^2}{3})$.
\item
$a=3^{6}.5^5.139$ and $b=2^{2}.3^5.5^5.139$. Let
$\t=\frac{\al}{15}$. Then $\t$ is integral with minimal
polynomial $F(X)=X^4+AX+B$ and $K=\q[\t]$,
 where $A=3^{3}.5^2.139$ and $B=2^{2}.3.5.139$. Thus
$2\z_K=(2,\t+1)(2,\t^2+\t+1)(2,\t)$, $3\z_K=(3,\t)^4$,
$5\z_K=(5,\t)^4$, $139\z_K=(139,\t)^4$. Since
$\triangle=2092367789117959822875=202317851.7.3^3.5^3.139^3.163$,
for every $p\not\in\{3,5,139\}$, $p\z_K$ is $p$-analogous to
$\bar F(X)$.
\end{enumerate}
\end{exs}

Lhoussain El Fadil  FPO, P.O. Box 638-Ouarzazte 45000, Morocco\\
lhouelfadil@hotmail.com
 \end{document}